\documentclass[12pt,psamsfonts]{amsart}
\usepackage{amsmath}
\usepackage{amssymb}
\usepackage[totalwidth=480pt,totalheight=680pt]{geometry}
\let\R\Real

\def\homeo#1#2#3
	{#1~\colon~#2~\stackrel\approxeq\to~#3}

\def\caption#1{\hfill \\
\hbox{}\hfil{\footnotesize #1}\hfil\hbox{}}
\def\ifnextchar#1#2#3{\let\tmpnce=#1%
    \def\tmpnca{#2}\def\tmpncb{#3}\futurelet\tmpncc\ifnch}%
\def\ifnch{\ifx\tmpncc\tmpnce\let\tmpncd=\tmpnca%
    \else\let\tmpncd=\tmpncb\fi\tmpncd}
\def\tpeinture #1 by #2 (#3){
  \vtop to #2{
    \hrule width #1 height 0pt depth 0pt
    \vfill
    \special{#3 y=#2}
    }
  }
\def\xpeinture #1 by #2 (#3){
  \hbox{$\vcenter to #2{
    \hrule width #1 height 0pt depth 0pt
    \vfill
    \special{#3 x=#1}
    }$}
  }
\def\ypeinture #1 by #2 (#3){
  \hbox{$\vcenter to #2{
    \hrule width #1 height 0pt depth 0pt
    \vfill
    \special{#3 y=#2}
    }$}
  }
\def\peinture{\ypeinture}
\def\bpeinture #1 by #2 (#3){
  \vbox to #2{
    \hrule width #1 height 0pt depth 0pt
    \vfill
    \special{#3 y=#2}
    }
  }
\def\Scaledpik[#1] #2 by #3 (#4){{ %
   \if#1t
       \tpeinture #2 by #3 (#4)%
   \else \if#1b
            \bpeinture #2 by #3 (#4)%
            \else \if#1x
                      \xpeinture #2 by #3 (#4)%
                      \else
                      \peinture #2 by #3 (#4)%
            \fi\fi
    \fi}}
\def\Scaledpij #1 by #2 (#3){{%
   \peinture #1 by #2 (#3) %
   }}
\def\scaledpicture{\ifnextchar[{\Scaledpik}{\Scaledpij}}
\def\centredpicture #1 by #2 (#3){
   \par\centerline{\hbox{
   \scaledpicture #1 by #2 (#3)}
   }}
\def\ifnextchar#1#2#3{\let\tmpnce=#1%
    \def\tmpnca{#2}\def\tmpncb{#3}\futurelet\tmpncc\ifnch}%
  \def\ifnch{\ifx\tmpncc\tmpnce\let\tmpncd=\tmpnca%
	\else\let\tmpncd=\tmpncb\fi\tmpncd}
   \def\pdiff#1#2{\dfrac{\partial#1}{\partial#2}}

\def\d{\operatorname{d}\!}

\def\modulus#1{\left|#1\right|}
\def\monthname{\ifcase\month\or Jan\or Feb\or March\or Apr\or %
    May\or June\or July\or Aug\or Sept\or Oct\or Nov\or Dec\fi}
\def\norm#1{\left\|#1\right\|}
      \def\normp#1_#2{\norm{#1}_{#2}}

\def\DEF{\stackrel{\scriptstyle\text{def}}{:=\!=}}

\input epsf.tex

\title[Analyzing biconcave vesicles]
{An analysis on the shape equation for biconcave axisymmetric vesicles}

\author[Thomas Au]{Thomas Kwok-keung Au}
\address{Department of Mathematics, The Chinese University of Hong Kong}
\email{thomasau@cuhk.edu.hk}

\author[Tom Wan]{Tom Yau-heng Wan}
\address{Department of Mathematics, The Chinese University of Hong Kong}
\email{tomwan@cuhk.edu.hk}

\newtheorem*{mthm}{Theorem}
\newtheorem{thm}{Theorem}[section]
\newtheorem{prop}[thm]{Proposition}
\newtheorem{lemma}[thm]{Lemma}
\newtheorem{cor}[thm]{Corollary}

\def\HelF{{\mathcal F}}
\def\bX{{\mathbf X}}

\def\tlambda{\lambda}
\def\tp{{p}}

\begin{document}
\begin{abstract}
We study the conditions on the physical parameters in the Helfrich
bending energy for lipid bilayer vesicles.  The variation equation
for embedded surface with a biconcave axisymmetric shape is analyzed
in detail.  This leads to simple conditions describing the 
solution and information about the geometry of the surface.
\end{abstract}
\maketitle
\setlength{\baselineskip}{28pt}

\section*{Introduction}

In this article, a vesicle is represented by a closed surface $\Sigma$ in
$\R^3$ with mean curvature $H$, surface area $\modulus{\Sigma}$ and it
encloses a volume $V$.  Its geometric shape is modelled by minimizing a
functional, sometimes called {\sl Helfrich functional},
$$
\HelF =
\int_\Sigma (2H+c_0)^2 \d S + \lambda\modulus{\Sigma} + p V
$$
with some physical constant parameters $c_0$, $\lambda$, and $p$.  The 
parameters carry the following meanings: $c_0$ is
the spontaneous curvature, $\lambda$ is the tensile stress, and
$p=p_o-p_i$ is the osmotic pressure difference between the outer
($p_o$) and inner ($p_i$) media.  We take a sign convention that
$H$ is negative for the standard sphere.

It has been observed experimentally long ago that a red blood cell is of
biconcave-discoid shape (which will be defined mathematically
later).  And it is a quest to find the appropriate
theoretical model for the bending energy.  Historically, the early model
of Canham, \cite{Canham}, is purely geometric and it is equivalent to the
Willmore functional, \cite[ch.~7]{Willmore}, which equals $\HelF$ with
$c_0 = \lambda = p = 0$, up to a constant.  Certainly, from a differential
geometer's view, this cannot be the correct model because the unique
minimum of the Willmore functional for topologically spherical vesicles is
the round sphere.  This is also observed by physicists
\cite{Helfrich1976}.  In fact, there are many important mathematical
studies of the Willmore functional because of its geometric implications;  
for instance, the existence of minimizers among a certain topological class
by Simon \cite{Simon}, and the conformal properties by Li and Yau
\cite{Li-Yau}.  We expect that their works may contribute to a certain
extent to a deeper theoretical understanding of the energy $\HelF$.

Helfrich takes physical condition together with the Gaussian
curvature into account and proposes a modified bending energy,
\cite{Helfrich1973}.  
The shape of blood cells and some other biological membranes
is closely related to the formation of lipid bilayer vesicles
in aqueous medium (e.g.~liquid crystal).
The physical condition is based on the elasticity of lipid bilayer 
vesicles.  According to the
Gauss-Bonnet Theorem, the integral of the Gaussian curvature is a
topological constant.  Thus, within a certain topological class of
$\Sigma$, Helfrich's bending energy can be reduced to $\HelF$ above.

Many properties of $\HelF$ are yet to be discovered, though
there are some experimental observations and numerical simulations,
\cite{Helfrich1977,Mutz-Bensimon}.  The existence and uniqueness of its
minimizer of a certain topology are still unknown.  It is also not known
whether the minimizer is symmetric in any sense.  Answers to these
questions require deep geometric analysis of the functional and studies in
this direction are rare.  Nevertheless, there are related works, such as,
on similar functionals, \cite{Simon,Li-Yau,Nitsche}; or on surface flows,
\cite{Elliott}.

The Euler-Lagrange equation corresponding to $\HelF$ is
$$
4\triangle_\Sigma H + 2(2H+c_0)\left[ 2(H^2-K)-c_0H \right] - 2\lambda H +
p = 0.
$$
In the past, much effort of physicists and biologists has been
spent on studying axisymmetric solution to this variational equation,
\cite{Helfrich1976,Luke,Seifert,OuYang-Helfrich1989,OuYang,Naito-Okuda-OuYang1,Naito-Okuda-OuYang2}.  
By an axisymmetric surface, we mean an embedded surface $\Sigma$ in $\R^3$
which is rotationally symmetric and has a reflection symmetry by the plane
perpendicular to the rotational axis.  With this additional assumption,
the fourth order equation can be reduced to a second order one, usually
referred to as the {\sl shape equation of axisymmetric vesicles}.  However, it
is still unknown to scientists how the solution depends on the
physical parameters and which parameters yield a solution corresponding to a
biconcave surface.  Our work is in this direction.  We can derive
conditions on the physical parameters for a solution having the
biconcave shape.  The conditions are easily expressed in terms of the
cubic and quadratic polynomials denoted by \begin{align*} Q(t) &= t^3 +
2c_0 t^2 + (c_0^2+\lambda)t - \frac{p}{2} ;\\ R(t) &= Q(t) - t^3.
\end{align*} Our main result can be stated as, \begin{mthm} For any $c_0$,
$\lambda$, and $p>0$ such that every real root of $Q$ is positive, there
exists axisymmetric biconcave surfaces which satisfies the Euler-Lagrange
equation of $\HelF$. \end{mthm} It should be remarked that when $c_0 > 0$,
$\lambda>0$, and $p>0$, the condition is always satisfied.  Helfrich's
numercial simulation produces a biconcave shape resembling a blood cell
when $c_0$ is positive.  Furthermore, we are also able to numerically
construct other interesting shapes when the condition is not satisfied.  
One is multiconcave and the other has no reflection symmetry.

This article is organized in the following way.  In~\S1, we first
give the differential equation for the revolving graph of
an axisymmetric solution.  The derivation of the equation is given
in the appendix.  Moreover, we formulate the problem of finding
special solution corresponding to a biconcave shape surface.  We also
present several variants of the equation which will be useful later.
The variants of the shape equation are
then studied in \S2 to show that our condition stated above is
sufficient for getting the expected special solution for the problem. 
The analysis and the estimates of geometric quantities are discussed
in detail here.  Finally, since the solution and its reflection no
longer form a graph at the reflection plane.  We will show that
the solution obtained in previous sections is still a solution across
the reflection symmetry.  It is sufficient to verify that it also
satisfies the variation equation at the reflection plane.
In the process, one more necessary geometric condition is obtained.

\subsubsection*{Acknowledgement}
We would like to thank our colleague
K.~S.~Chou, who has been encouraging in our project and making valuable
suggestions.

\setcounter{section}{0}

\section{The Equation, the Problem, and the Conditions}
An axisymmetric surface is a closed embedded surface $\Sigma$ in
$\R^3$ with a rotation symmetry and a reflection symmetry by the
plane perpendicular to the rotation axis.  It is biconcave if there
are exactly two components of negative Gaussian curvature.
Without loss of generality, the rotational axis is labelled the $z$-axis
and the plane of reflection is the $xy$-plane.  Then the surface
$\Sigma$ can be obtained by revolving a radial curve about the $z$-axis
on the upper half plane and reflecting it to the lower half.
Typically, a biconcave one is obtained by revolving and reflecting a
curve shown in the picture below.
\begin{center}
\mbox{\epsfxsize=8cm \epsfbox{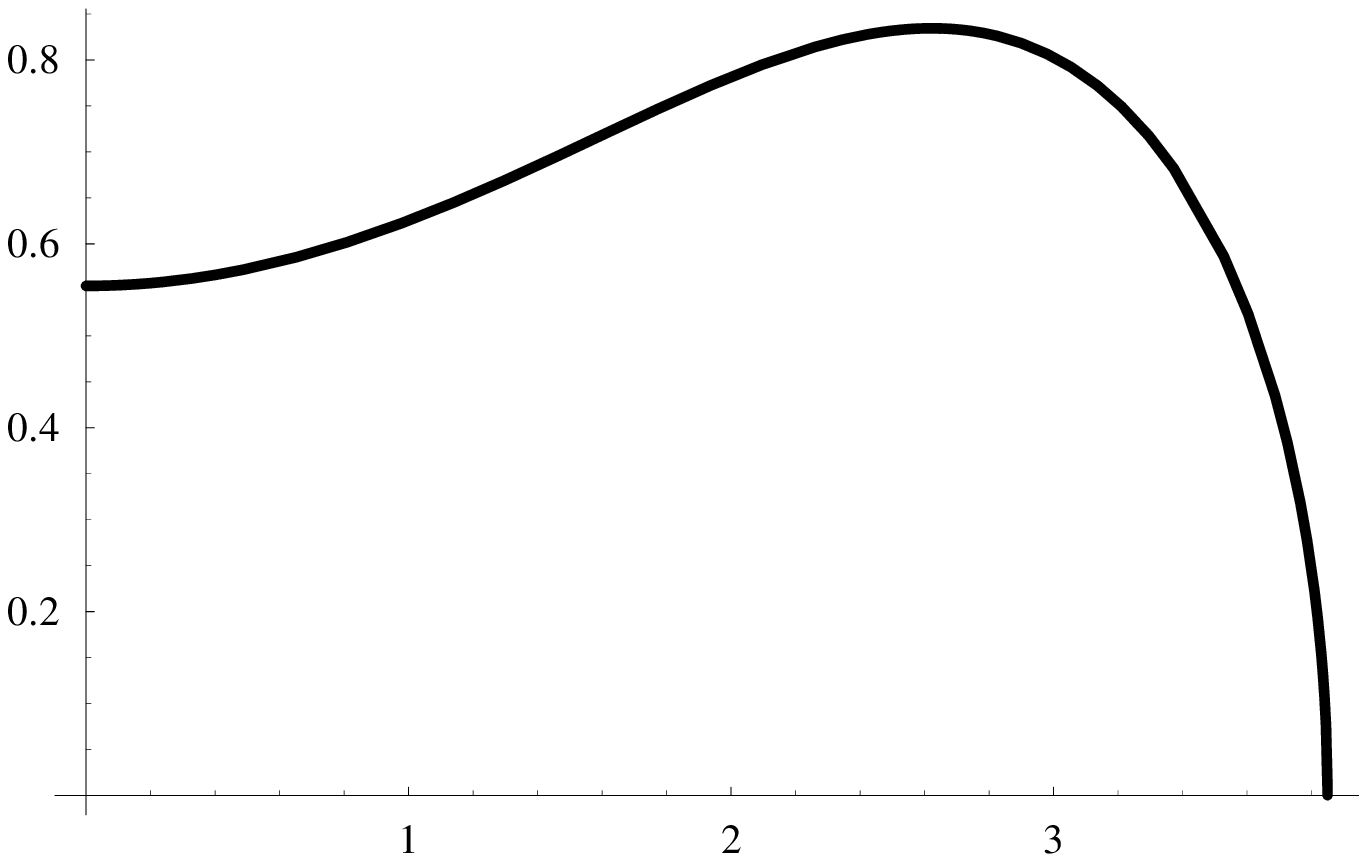}} \vspace*{-4ex} \\
{\footnotesize A cross-section of a biconcave axisymmetric surface
\vspace*{-2ex} \\
(with $c_0=1$, $\tlambda=0.25$, $\tp=1$)}
\end{center}

In other words, one may parametrize the upper part of $\Sigma$ by
$\bX = (r\cos\theta, r\sin\theta, z(r)),$
with a function $z(r)$ defined for $r$ in some interval
$[0,r_\infty]$, where $z(r)>0$ for $r\in [0,r_\infty)$ and $z(r_\infty)=0$.
There are natural boundary conditions imposed by the rotation and
reflection symmetries of the surface.  The obvious ones are $z'(0)=0$
and $z'(r)\to -\infty$ as $r\to r_\infty$.  However, more subtle ones
arise from the regularity at the end-point $r_\infty$, which we will
discuss later in \S 3.
For a biconcave surface, there is $r_M < r_\infty$ such that $z''(r) > 0$
if and only if $r \in [0,r_M)$.  This is equivalent to require that $z'$
has only a unique maximum and no other critical point.

For a biconcave axisymmetric surface, the variation equation for the
Helfrich functional is reduced to a second order ordinary differential
equation, traditionally called the shape equation.
With the notation $w(r) = z'(r)$, the shape equation is
\begin{equation} \label{VEqn1}
\begin{aligned}
\frac{2r}{(1+w^2)^{5/2}}{w''} =  &\frac{5rw}{(1+w^2)^{7/2}}{w'}^2 - \frac{2{w'}}{(1+w^2)^{5/2}} \\
&\qquad {} + \frac{2w+w^3}{r(1+w^2)^{3/2}} + \frac{2c_0 w^2}{1+w^2} + \frac{(c_0^2+\tlambda)r w}{(1+w^2)^{1/2}} - \frac{\tp r^2}{2}.
\end{aligned}
\end{equation}
We leave the derivation of this equation to the appendix to focus on
the idea and analysis of the equation.

We are going to study special solution $w(r)$ to the initial value
problem on this variational equation
with initial choice $w(0)=0$ and $w'(0)={w_0'} > 0$.  Specifically, we
look for a solution with the additional requirements that
\begin{description}
\item[C1]
it is unimodal; and
\item[C2]
there is a finite number $r_\infty$
such that $w(r)\to -\infty$ as $r\to r_\infty$; and
\item[C3]
$\displaystyle -\infty < \int_0^{r_\infty} w(r)\d r < 0.$
\end{description}

In \cite{Helfrich1976,Julicher-Seifert,Zheng-Liu}, there are versions of
the same equation written in term of the angle $\psi$ between the
surface tangent and the plane perpendicular to rotational axis.  However,
it is convenient for our discussion to write the equation in the above form.

With our notation of the polynomials $Q$ and $R$, after multiplying with
$r {w'}$, the equation is equivalent to,
\begin{align}
\left[ \frac{r^2 {w'}^2}{(1+w^2)^{5/2}} \right]' &= \left[ \frac{w^2}{(1+w^2)^{1/2}} \right]' + r^3 {w'} R(\kappa(r)); \quad\mbox{or}\tag{\ref{VEqn1}a} \label{VEqn2} \\
\left[ \frac{r^2 {w'}^2}{(1+w^2)^{5/2}} \right]' &= \left[ \frac{-2}{\sqrt{1+w^2}} \right]' + r^3 {w'} Q(\kappa(r)), \tag{\ref{VEqn1}b} \label{VEqn3}
\end{align}
where $\kappa(r) = \dfrac{w}{r\sqrt{1+w^2}}$.  It will be seen that these
groupings of the lower order terms are important in the analysis the equation.

The understanding of $\kappa(r)$ also provides useful information about
the solution.  First, its derivatives are given by,
\begin{align*}
\kappa'(r) &= \frac{w'}{r(1+w^2)^{3/2}} - \frac{w}{r^2\sqrt{1+w^2}}, \\
\kappa''(r) &= \frac{w''}{r(1+w^2)^{3/2}} - \frac{3w{w'}^2}{r(1+w^2)^{5/2}} - \frac{2w'}{r^2(1+w^2)^{3/2}} + \frac{2w}{r^3\sqrt{1+w^2}}.
\end{align*}
Then we have the equation for $\kappa(r)$,
\begin{equation} \label{VEqnk}
r\kappa'' = \frac{-r\kappa(r\kappa'+\kappa)^2}{2(1-r^2\kappa^2)} - 3\kappa' + \frac{r Q(\kappa(r))}{2(1-r^2\kappa^2)}.
\end{equation}

Let us end this section by remarking on a few geometric
quantities in terms of $w$.  Firstly, $\kappa = \dfrac{w}{r\sqrt{1+w^2}}$
is the principle curvature in the meridinal (rotational) direction. 
Another principle curvature, the longitudinal one, is given by
$\dfrac{w'}{(1+w^2)^{3/2}}$, which occurs in $\kappa'$.  Their product is
the Gaussian curvature, which is expected to be positive at the rotational
axis and reflection plane but negative somewhere along the circle
defined by the zero of $w$.

\section{Analysis of the Equation}
In this section, we prove the conditions for the existence of the required
special solution to the initial
value problem on equation~(\ref{VEqn1}) described in the preceding
section.  The strategy is the following analysis on the equation.  We
study the principle curvature $\kappa$, which is positive initially (at the
rotational axis), i.e., $\kappa(0)=w_0'>0$.  If $w_0'$ is not too large in
terms of the roots of $Q$, $\kappa$ must decrease and eventually becomes zero
at $r_0$ for some $0<r_0<r_\infty$.  After $w$ becomes negative, it
continues its descent and blows down to $-\infty$ at finite distance
$r_\infty$.  This guarantees {\bf C1} that $w$ is unimodal and has a unique
zero at $r_0$.  In order to verify that the solution satisfies the
requirements {\bf C2} and {\bf C3}, we establish, in terms of
$w_0'$, the estimates on $r_0$, $r_\infty$ and values of $w$, $w'$ at
these positions.

\begin{lemma}\label{lem-kdrops}
Let ${w_0'} > 0$ and $R(t) < 0$ for $t\in[0,w_0']$.  If $w \geq 0$ on an interval $[0,r_0)$, then $\kappa(r)$ decreases on $(0,r_0)$.
\end{lemma}
\begin{proof}
By continuity of $\kappa$ and that $\displaystyle \lim_{r\to 0}\kappa(r) = {w_0'}>0$, there exists $\varepsilon>0$ such that $R(\kappa(r)) < 0$, $w \ge 0$, and ${w'} > 0$ on $(0,\varepsilon)$.
Integration of the equation (\ref{VEqn2}) on $(0,\varepsilon)$ gives
$$
\frac{r^2 {w'}^2}{(1+w^2)^{5/2}} < \frac{w^2}{\sqrt{1+w^2}}
\quad \mbox{on}\quad (0,\varepsilon).
$$
Thus, ${\kappa'}(r) = \dfrac{{w'}}{r(1+w^2)^{3/2}} - \dfrac{w}{r^2\sqrt{1+w^2}} < 0$ for $r\in (0,\varepsilon)$. 
This also implies that $\kappa(r)$ remains in the interval $(0,w_0']$ for $r\in (0,\varepsilon)$ and hence the argument works before we hit the first zero $r_M$ of $w'$, i.e., the first critical point of $w$.
This shows that $r_1=\inf\{ r\in(0,r_0)\, :\, \kappa'(r)\ge 0 \}\ge r_M>0$ if the set is nonempty.

Now by the smoothness of $\kappa$, we conclude that $\kappa(r_1)\in[0,w_0']$, $\kappa'(r_1)=0$ and $\kappa''(r_1)\ge 0$. Putting this into the equation (\ref{VEqnk}), we have
$$
0\le r_1\kappa''(r_1)=\frac{r_1 R(\kappa(r_1))}{2\left( 1- r_1^2\kappa^2(r_1) \right)} <0.
$$
This is a contradiction and hence the set $\{ r\in(0,r_0)\, :\, \kappa'(r)\ge 0 \}$ must be empty.
\end{proof}

By the same token, we can also establish a criterion for a growing solution. Since we don't need it for our further discussion, we will omit the proof.
\begin{prop}
If $c_0>0$, $p>0$, ${w_0'} > 0$ and $R({w_0'}) > 0$, then $w$ is increasing  and blows-up to $+\infty$.
\end{prop}

For our future discussion, we denote a few quantities which depends only on the polynomial $Q$ and the initial data $w_0'$ as follow:
\begin{align*}
\mu (w_0')&= -\min\left\{ Q(t)~:~ 0 \leq t \leq w_0' \right\} 
= \max\left\{ -Q(t)~:~ 0 \leq t \leq w_0' \right\} ;\\
\delta_{+}(w_0') &= \min\left\{ -Q(t)~:~ 0 \leq t \leq w_0' \right\} = -\max\left\{ Q(t)~:~ 0 \leq t \leq w_0' \right\} ;\\
\delta_{-} &= \min\left\{ -Q(t)~:~ t \leq 0 \right\}
= -\max \left\{ Q(t)~:~ t \leq 0 \right\}.
\end{align*}
Note that $\displaystyle\lim_{w_0'\to 0}\delta_{+}(w_0')=p/2=\lim_{w_0'\to 0}\mu (w_0')$, and $\delta_{-}$ is independent of $w_0'$.

We will now show that under reasonable condition, the solution will not blow up to $+\infty$.  It is because the cubic lower terms are dominated by $\delta_{+}$.
\begin{lemma} \label{thm-pos}
Suppose that all real roots of $Q$ are positive.  
If $Q<0$ on $[0,{w_0'}]$, i.e., $w_0'$ less than the smallest real root of $Q$, then there is $r_0 > 0$ with $r_0^2 < \frac{16{w_0'}}{\delta_+}$ such that on the interval $(0,r_0)$, 
$$
w>0,\quad w(r_0) = 0,\quad  \mbox{and}\quad w'(r_0)<-\frac{\delta_{+}}{8}r_0^2.
$$
In addition, if $64 {w_0'}^3<27\delta_+$, then 
$$
\int_0^{r_0} w(r) \d r \le \frac{4{w_0'}^2}{\delta_{+}^{3/2}}
\left(1-\frac{64 {w_0'}^3}{27\delta_+} \right)^{-1/2}.
$$
\end{lemma}
\begin{proof}
Consider the principal curvature $\kappa(r) = \dfrac{w}{r\sqrt{1+w^2}}$ and rewrite equation~(\ref{VEqn1}) or (\ref{VEqn2}) as
$$
{\kappa''} = \frac{-3{\kappa'}}{r} - \frac{w{w'}^2}{r(1+w^2)^{5/2}} + \frac{1+w^2}{2}Q(\kappa(r)).
$$
By lemma~\ref{lem-kdrops}, ${\kappa'}(r) < 0$ on a neighborhood of $0$ and hence $\kappa(r)<w_0'$ near $0$.  Using the assumption on $Q$ and the positivity of $w$ near $0$, we have
$$
r^3{\kappa''} + 3r^2{\kappa'} \leq \frac{-\delta_+}{2}r^3.
$$
Therefore, for small $r$,
\begin{equation}\label{ineq-upkprime}
\kappa'(r)\le -\frac{\delta_+}{8}r^2,
\end{equation}
and
\begin{equation}\label{ineq-upk}
\kappa(r) \leq {w_0'} - \frac{\delta_+}{16}r^2.
\end{equation}
These inequalities show that $\kappa$ remains in $(0, w_0')$ as long as $w\ge 0$.  In turns, they themselves hold as long as $w \geq 0$.  Clearly, they guarantee the existence of $r_0$ with $r_0^2\leq \dfrac{16 {w_0'}}{\delta_+}$ such that $\kappa(r_0) = 0$ and so $w(r_0)=0$. Then by (\ref{ineq-upkprime}), we have $w'(r_0)<-\delta_{+}r_0^2/8$. This proves the first statement of the theorem.

To prove the second statement, we temporarily let $a = {w_0'}$ and $b=\delta_+/16$.  Then on $[0,r_0)$,
$$
\frac{w}{r\sqrt{1+w^2}} \leq a - b r^2.
$$
With the assumption that $4a^3 < 27 b$, it can be written as
$$
w \leq \frac{r(a-br^2)}{\sqrt{1-r^2(a-br^2)^2}}.
$$
Therefore,
$$
w \leq \frac{r(a-br^2)}{\sqrt{1-\frac{4a^3}{27b}}},
$$
and
\begin{align*}
\int_0^{r_0} w \d r & \leq \int_0^{\sqrt{a/b}} \frac{r(a-br^2)}{\sqrt{1-\frac{4a^3}{27b}}}\d r  \\
&\leq  \frac{1}{\sqrt{b-\frac{4a^3}{27}}} \frac{a^2}{4b} 
= \frac{4{w_0'}^2}{\delta_{+}} \left(1 - \frac{64{w_0'}^3}{27\delta_+}\right)^{-1/2}.
\end{align*}
This completes the proof of the lemma.  
\end{proof}
Immediately from the lemma we have the following
\begin{cor}
Suppose that all real roots of $Q$ are positive. Then
$$
\limsup_{w_0'\to 0}\frac{r_0^2}{w_0'} \le \frac{32}{p}\quad\mbox{and}\quad
\limsup_{w_0'\to 0}\frac{1}{{w_0'}^2}\int_0^{r_0}w\d r \le \frac{8}{p}.
$$
\end{cor}
In the above, we obtained an upper bound for $r_0$, where $w$ first hits zero.  Next, a lower bound is established, which is essential for future estimates. We note that from lemma \ref{thm-pos}, $w$ has at least one maximum in the interval $(0,r_0)$. Since $w$ is increasing at the beginning, the first critical point must be a maximum point.  From now on, we let $\xi = \xi(w_0') = 1 - \dfrac{64{w_0'}^3}{27\delta_{+}} > 0$.  
\begin{lemma} \label{lem-lowr0}
Let $r_M\in (0,r_0)$ be the first critical point of $w$.  Then
$$
\liminf_{w_0'\to 0}\frac{r_0^2}{w_0'} \ge 
\liminf_{w_0'\to 0}\frac{r_M^2}{w_0'} \geq \frac{32}{3p}, \qquad\qquad
\liminf_{w_0'\to 0}\frac{w'(r_0)}{w_0'} \ge -2.
$$
\end{lemma}
\begin{proof}
Using equation (\ref{ineq-upk}), $\kappa >0$ on $[0, r_0)$, and $r_0^2< 16w_0'/\delta_+$, one can check that
$$
1 - r^2\kappa(r)^2 \geq \xi \quad\mbox{on } (0,r_0).
$$
Putting together the assumption on $Q$, the fact that $\kappa<w_0'$ on $(0, r_0)$, and the inequality~(\ref{ineq-upkprime}) into the equations (\ref{VEqnk}), one obtains
\begin{align*}
\left( r\kappa\right)'' 
&\geq \frac{-1}{2\xi}r\kappa(r\kappa'+\kappa)^2 + \frac{1}{8}\delta_{+}r - \frac{\mu}{2\xi}r.\\
&\geq \frac{-1}{2\xi}r{w_0'}\left[(r\kappa)'\right]^2 + \frac{1}{8}\delta_{+}r - \frac{\mu}{2\xi}r \\
&= \frac{-w_0'}{2\xi} \left[ {(r\kappa)'}^2 + \frac{\mu - \frac{\xi}{4}\delta_{+}}{{w_0'}} \right]r.
\end{align*}
Note that in the above,
$\mu = \max\left\{-Q(t): 0\leq t \leq w_0' \right\} \geq \delta_{+}$.

Let $y = (r\kappa)' = \dfrac{w'}{(1+w^2)^{3/2}}$ and $A^2 = \dfrac{\mu - \frac{\xi}{4}\delta_{+}}{{w_0'}} \geq \dfrac{\left(1 - \frac{\xi}{4}\right)\delta_{+}}{{w_0'}} > \dfrac{3\delta_{+}}{4w_0'} > 0$.  We then have
$$
\frac{y'}{y^2 + A^2} \geq \frac{-w_0'}{2\xi}r.
$$
Integrating from $0$, it yields
\begin{equation}\label{ineq-lowr0}
\arctan \left(\frac{y}{A}\right) \geq \arctan \left(\frac{{w_0'}}{A}\right) - \frac{Aw_0'}{4\xi}r^2.
\end{equation}
Taking $r=r_M$ and multiplying by $A/w_0'$, as $y(r_M)=0$ and $r_0 \geq r_M$, we have
\begin{align*}
A^2 r_0^2 \geq A^2 r_M^2 &\geq 4\xi\cdot\frac{A}{{w_0'}}\arctan \left(\frac{{w_0'}}{A}\right) \geq 4\xi\left(1-\frac{{w_0'}^2}{3A^2}\right). \\
\intertext{Then,}
\frac{r_0^2}{w_0'} \geq \frac{r_M^2}{w_0'} &\geq \frac{4\xi} {\left(\mu - \frac{\xi}{4}\delta_+\right)} \left[ 1 - \frac{{w_0'}^3} {3\left(\mu - \frac{\xi}{4}\delta_+\right)} \right].
\end{align*}
Considering the limiting situation, by $\displaystyle\lim_{w_0'\to 0}\mu=\lim_{w_0'\to 0}\delta_+=p/2$ and $\displaystyle\lim_{w_0'\to 0}\xi = 1$, the first result follows.

For the second estimate, we first take $r\to r_0$ in (\ref{ineq-lowr0}), apply $r_0^2 < 16w_0'/\delta_+$ to obtain
\begin{align*}
\frac{A(w'(r_0)-w_0')}{A^2 + w'(r_0)w_0'} &\geq \tan\left(\frac{-4A{w_0'}^2}{\delta_+\xi}\right) \geq \frac{-4A{w_0'}^2}{\delta_+\xi} \cdot \left[ \cos\left(\frac{4A{w_0'}^2}{\delta_+\xi}\right) \right]^{-1}. \\
\intertext{Then,}
\frac{w'(r_0)}{w_0'} &\ge 
\frac{ \cos\left(\frac{4A{w_0'}^2}{\delta_+\xi}\right) - \frac{4A^2{w_0'}}{\delta_+\xi}} { \cos\left(\frac{4A{w_0'}^2}{\delta_+\xi}\right) +\frac{4{w_0'}^3}{\delta_+\xi}}.
\end{align*}
This gives the second estimate immediately by letting $w_0'\to 0$.
\end{proof}
Combining the previous lemmas, we have the following
\begin{cor}\label{cor-bothr0}
For any $\varepsilon >0$, there exists $\beta>0$ such that if $w_0'\in (0,\beta)$, then
$$
\frac{32}{3p}-\varepsilon \le \frac{r_M^2}{w_0'} \le \frac{r_0^2}{w_0'}
\le \frac{32}{p}+\varepsilon
$$
and
$$
-2-\varepsilon \le \frac{w'(r_0)}{w_0'} \le -\frac{2}{3}+\varepsilon
$$
\end{cor}
The $r_M$ that we discussed in the previous lemma is in fact the unique critical point of $w$ in $(0, r_0)$ provided $w_0'$ is sufficiently small.
\begin{lemma}
Under conditions of lemmas~{\em\ref{thm-pos}} and {\em\ref{lem-lowr0}}, if $w_0'$ is sufficiently small, then $w$ has a unique maximum at $r_M \in [0,r_0)$ and has no other critical point.
\end{lemma}
\begin{proof} Suppose $w$ has another critical points in $(0,r_0)$, then according to the definition of $r_M$, $w$ attains a positive nonmaximum critical at $r_m\in (r_M, r_0)$, i.e., $w'(r_m)=0$ and $w''(r_m)\geq 0$.
By the corollary~\ref{cor-bothr0}, given any $\varepsilon >0$, $r_m^2 > r_M^2 \ge \left(\frac{32}{3p}-\varepsilon\right)w_0'$ for sufficiently small $w_0'$. 
On the other hand, substituting $w'(r_m)=0$ and $w''(r_m)\geq 0$ into the equation~(\ref{VEqn1}), we have
\begin{align*}
0 &\leq \frac{2w_m}{r_m^3(1+{w_m}^2)^{3/2}} + Q(\kappa(r_m)) \\
&\leq \frac{2}{r_m^2}\kappa(r_m) + Q(\kappa(r_m)),
\end{align*}
where $w_m=w(r_m)$.
Since $0<\kappa(r_m) \leq w_0'$, thus, for sufficiently small $w_0'$, we also have $Q(\kappa(r_m)) < -p/2 + \varepsilon$.  This leads to
$$
0 \leq \frac{2}{\frac{32}{3p}-\varepsilon} - \frac{p}{2} + \varepsilon
$$
which is clearly a contradiction since $\varepsilon$ is arbitrary.
\end{proof}

Next, we claim that under our condition on $Q$, after hitting $0$, $w$ decreases and goes to negative infinity in finite distance.
\begin{thm} \label{thm-neg}
Suppose that all real roots of $Q(t)$ are positive and $r_0 > 0$ is given by lemma~{\em\ref{thm-pos}}.  Then ${w'} < 0$ for $r \geq r_0$ and there is a finite number $r_\infty > r_0$ such that
$\displaystyle \lim_{r\to r_\infty} w(r) = -\infty$.  Moreover, we have
$$
r_\infty - r_0 \leq \frac{\pi}{2\sqrt{\delta r_0^2 \left|w'(r_0)\right|}}, $$
where $\delta = \min\left\{\frac{\delta_{+}}{8}, \frac{\delta_{-}}{2}\right\}$ and $\delta_{-}=\min\left\{ -Q(t)~:~ t \leq 0 \right\}$.
\end{thm}
\begin{proof}
For convenience, let us temporarily denote $v = -w$ for $r > r_0$.  Since $\kappa(r) < 0$ and ${v'} > 0$ on $(r_0, r_0+\varepsilon)$ for some $\varepsilon>0$, we have $-{v'} Q(\kappa(r)) > \delta_{-}{v'}$ on the same interval and the equation~(\ref{VEqn3}) gives
\begin{equation} \label{VEqnneg}
\left[ \frac{r^2{v'}^2}{(1+v^2)^{5/2}} \right]' = \left[ \frac{-2}{(1+v^2)^{1/2}} \right]' - r^3 {v'} Q(\kappa(r)) \geq \left[ \frac{-2}{(1+v^2)^{1/2}} \right]' + \delta_{-} r^3 {v'}
\end{equation}
on $(r_0, r_0+\varepsilon)$.
In particular,
$$
\left[ \frac{r^2{v'}^2}{(1+v^2)^{5/2}} \right]' \geq \left[ \frac{-2}{(1+v^2)^{1/2}} \right]',
$$
which gives
$$
\frac{r^2{v'}^2}{(1+v^2)^{5/2}} - r_0^2 {v'}(r_0)^2 \geq \frac{-2}{(1+v^2)^{1/2}} + 2 \geq 0 .
$$
This implies that $v'$ does not vanish and therefore
\begin{equation}\label{ineq-lowvprime}
r{v'} \geq \frac{r{v'}}{(1+v^2)^{5/4}} \geq r_0 {v'}(r_0) > 0.
\end{equation}
This, in turns, implies that $v$ is increasing, $\kappa$ remains negative and (\ref{VEqnneg}) holds as long as $v$ is defined.
Substitute (\ref{ineq-lowvprime}) into the equation (\ref{VEqnneg}) again, it gives
$$
\left[ \frac{r^2{v'}^2}{(1+v^2)^{5/2}} \right]' \geq  \delta_{-} r^3 {v'} \geq \delta_{-}r_0 {v'}(r_0)r^2 \geq \delta_{-}r_0^2 {v'}(r_0)r
$$
for all $r>r_0$ such that $v$ is defined. 
Now, integrating the above from $r_0$, we conclude that
$$
\frac{r^2{v'}^2}{(1+v^2)^{5/2}} \geq r_0^2v'(r_0)^2 + \frac{\delta_{-}r_0^2 {v'}(r_0)}{2}(r^2-r_0^2).
$$
In the case that $\delta_{+} \geq 4\delta_{-}$, we have $v'(r_0)\ge \delta_{-}r_0^2/2$ by lemma~\ref{thm-pos}
\begin{align*}
\frac{r^2{v'}^2}{(1+v^2)^{5/2}} &\geq r_0^2v'(r_0)\left(v'(r_0) - \frac{\delta_{-}}{2}r_0^2\right) + \frac{\delta_{-}r_0^2 {v'}(r_0)}{2}r^2 \\
&\geq \frac{\delta_{-}r_0^2 {v'}(r_0)}{2}r^2.
\end{align*}
Otherwise, for $\zeta = \frac{\delta_{+}}{4\delta_{-}}<1$, one has $r_0^2v'(r_0)^2 \geq \frac{1}{2}\zeta\delta_{-}v'(r_0)r_0^4.$  It follows that
\begin{align*}
\frac{r^2{v'}^2}{(1+v^2)^{5/2}} &\geq r_0^2v'(r_0)^2 + \frac{\delta_{-}r_0^2 {v'}(r_0)}{2}(r^2-r_0^2) \\
&= r_0^2v'(r_0)^2 + \frac{\zeta \delta_{-}r_0^2 {v'}(r_0)}{2}(r^2-r_0^2) + \frac{(1-\zeta) \delta_{-}r_0^2 {v'}(r_0)}{2}(r^2-r_0^2) \\
&\geq r_0^2v'(r_0)\left[v'(r_0) - \frac{\zeta \delta_{-}{v'}(r_0)}{2}r_0^2\right] + \frac{\zeta \delta_{-}r_0^2 {v'}(r_0)}{2}r^2. \\
\intertext{Therefore,}
\frac{r^2{v'}^2}{(1+v^2)^2} &\geq \frac{r^2{v'}^2}{(1+v^2)^{5/2}} \geq \frac{\zeta \delta_{-}r_0^2 {v'}(r_0)}{2}r^2 = \frac{1}{8}\delta_{+}r_0^2v'(r_0)r^2.
\end{align*}
Hence, with $\delta = \min\left\{\frac{1}{4}\delta_{+},
\frac{1}{2}\delta_{-}\right\},$ one has
$$
\frac{ {v'}}{1+v^2}\ge \frac{ {v'}}{(1+v^2)^{5/4}} \geq  \sqrt{\delta r_0^2{v'}(r_0)}
$$
as long as $v$ is defined.
Further integrating from $r_0$ gives
$$
\arctan v \geq  \sqrt{\delta r_0^2{v'}(r_0)} (r-r_0).
$$
This clearly shows that $v$ blows up before $r$ goes to infinity.  Furthermore, the upper bound of $r_\infty$ follows.
\end{proof}
\begin{cor}\label{cor-uprinf}
Under the same conditions of theorem~{\em\ref{thm-neg}}, we have,
$$
\infty > \int_{r_0}^{r_\infty} \modulus{w(r)} \d r  \geq B^{-1}\log\left( \frac{1}{\cos(B(r_\infty-r_0))}\right) \geq \frac{B}{2}(r_\infty-r_0)^2,
$$
where $B = \sqrt{\delta r_0^2 v'(r_0)} \geq \sqrt{\frac{\delta\delta_{+}}{8}}r_0^2 \geq \frac{1}{\sqrt{2}}\delta r_0^2$.
\end{cor}
\begin{proof}
From (\ref{VEqnneg}), by simply dropping the term $-r^3v'Q(\kappa(r)) > 0$, we have
$$
\frac{r^2{v'}^2}{(1+v^2)^{5/2}} \geq 2\left(1-\frac{1}{\sqrt{1+v^2}}\right).
$$
After multiplying by $\sqrt{1+v^2}(\sqrt{1+v^2}+1)$, it becomes
\begin{align*}
(\sqrt{1+v^2}+1)\frac{r^2{v'}^2}{(1+v^2)^2} & \geq 2(\sqrt{1+v^2}+1)(\sqrt{1+v^2}-1) = 2v^2. \\
\intertext{Moreover, $2\sqrt{1+v^2} > (\sqrt{1+v^2}+1)$, so one has}
r_\infty^2 \cdot \frac{{v'}^2}{(1+v^2)^{3/2}} &\geq v^2. \\
\intertext{Taking square root and integrating, it yields}
r_\infty\int_{r_0}^{r_\infty} \frac{v'}{(1+v^2)^{3/4}} \d r &\geq \int_{r_0}^{r_\infty} v \d r ,
\end{align*}
in which the left hand side is obviously convergent.  The lower bound follows easily from the theorem.
\end{proof}

It can be seen that the area of negative part has a lower bound of order $w_0'(r_\infty-r_0)^2$.  We wish to establish that this dominates the area of the positive part.  Let us consider the equation of $v(r)=-w(r)$ for $r > r_0$ again, i.e.,
\begin{align*}
\left[ \frac{r^2 {v'}^2}{(1+v^2)^{5/2}} \right]' &= \left[ \frac{v^2}{(1+v^2)^{1/2}} \right]' - \frac{2c_0 v^2v' r}{1+v^2} + \frac{(c_0^2+\lambda)r^2 vv'}{(1+v^2)^{1/2}} + \frac{pr^3v'}{2} \\
&\leq \left[ \frac{v^2}{(1+v^2)^{1/2}} \right]' + \left[ \frac{2|c_0|v^2r_\infty}{1+v^2} + \frac{(c_0^2+\lambda)r_\infty^2 v}{(1+v^2)^{1/2}} + \frac{pr_\infty^3}{2} \right]v' \\
&\leq \left[ \frac{v^2}{(1+v^2)^{1/2}} \right]' + \left[ 2|c_0|r_\infty + (c_0^2+\lambda)r_\infty^2 + \frac{pr_\infty^3}{2} \right]v'.
\end{align*}
Therefore, after integrating from $r_0$,
\begin{align*}
\frac{r^2 {v'}^2}{(1+v^2)^{5/2}} &\leq r_0^2v'(r_0)^2 + \frac{v^2}{(1+v^2)^{1/2}} + \left[2|c_0|r_\infty + (c_0^2+\lambda)r_\infty^2  + \frac{pr_\infty^3}{2}\right]v \\
&\leq r_0^2v'(r_0)^2 + \left[1 + 2|c_0|r_\infty + (c_0^2+\lambda)r_\infty^2 + \frac{pr_\infty^3}{2}\right]v \\
\frac{r_0^2 {v'}^2}{(1+v^2)^{5/2}}&\leq  r_0^2v'(r_0)^2 + \left[ 1 + 2|c_0|r_\infty + (c_0^2+\lambda)r_\infty^2 + \frac{p}{2}r_\infty^3\right]v.
\end{align*}
From this, we have proved a comparison between the growth of $w'$ and $\modulus{w}$.
\begin{prop}\label{prop-wprimeord}
If $w(r)$ is a solution to equation~{\em (\ref{VEqn1})} which blows down to $-\infty$ at $r_\infty$, then
$
\dfrac{{w'}^2}{\modulus{w}(1+w^2)^{5/2}}
$
is bounded.
\end{prop}
Continue with the analysis, after integration and applying Holder inequality, we have
\begin{align*}
r_0^2 &\leq r_0^2\left(\int_{r_0}^{r_\infty} \frac{v'}{(1+v^2)^{5/4}}\cdot 1\right)^2 \leq r_0^2 \int_{r_0}^{r_\infty} \frac{{v'}^2}{(1+v^2)^{5/2}} \int_{r_0}^{r_\infty}1  \\
&\leq  r_0^2v'(r_0)^2 (r_\infty -r_0)^2 + (r_\infty -r_0)\left[ 1 + 2|c_0|r_\infty + (c_0^2+\lambda)r_\infty^2 + \frac{p}{2}r_\infty^3\right] \int_{r_0}^{r_\infty} v.
\end{align*}
Letting $x=r_\infty -r_0$ and using \ref{cor-bothr0} that $r_0^2<\frac{64}{\delta_+}w_0'$ and $v'(r_0)<2w_0'$, we have
\begin{equation}\label{ineq-Intv}
r_0^2 \leq  \frac{64}{\delta_+}{w_0'}^3 x^2 + x\left[ a_0+a_1x+a_2x^2 + \frac{p}{2}x^3\right] \int_{r_0}^{r_\infty} v,
\end{equation}
where
\begin{align*}
a_0&=1+2|c_0|r_0+(c_0^2+\lambda)r_0^2+\frac{p}{2}r_0^3,\\
a_1&=2|c_0|+2(c_0^2+\lambda)r_0+\frac{3p}{2}r_0^2 \mbox{ and}\\
a_2&=(c_0^2+\lambda)+\frac{3p}{2}r_0.
\end{align*}
We now claim that there exists $C>0$ independent of $w_0'$ such that $\int_{r_0}^{r_\infty} v\ge Cw_0'$ for  sufficiently small $w_0'$. Otherwise, for any $\varepsilon>0$, there is a $w_0' <\varepsilon$ such that
$$
\varepsilon w_0' > \int_{r_0}^{r_\infty} v.
$$
Then by corollary~\ref{cor-uprinf}, we have
$$
\varepsilon w_0' >\frac{\delta}{2\sqrt{2}}r_0^2x^2.
$$
Using the lower estimate of $r_0^2$ in terms of $w_0'$, we conclude that
$$
x^2 < \varepsilon \frac{2\sqrt{2}}{\delta}\frac{w_0'}{r_0^2} \le 
\frac{2\sqrt{2}\delta_+}{\delta}\varepsilon.
$$
Putting this into (\ref{ineq-Intv}), we have
\begin{align*}
c_1w_0' &\le c_2{w_0'}^3\varepsilon +c_3\varepsilon^{1/2} \left\{
1+o({w_0'}^{1/2})+\left[2|c_0|+o({w_0'}^{1/2})\right]c_3\varepsilon^{1/2} + \right. \\
& \qquad \left. {} + \left[ c_0^2+\lambda + o({w_0'}^{1/2}) \right]c_3^2\varepsilon +\frac{p}{2}c_3^3\varepsilon^{3/2}
 \right\}\varepsilon w_0',
\end{align*}
for some positive constants $c_1$, $c_2$ and $c_3$. This is impossible since this implies
$$
0<c_1\le c_2{w_0'}^2\varepsilon +c_3\varepsilon^{3/2}\left[ 1 + o(\varepsilon^{1/2}) + o({w_0'}^{1/2}) \right] \to 0,
$$
which is a contradiction.

\begin{thm}
With the conditions of the theorems {\em\ref{thm-pos}, \ref{lem-lowr0}, and \ref{thm-neg}}, we have
$$
-\infty < \int_0^{r_\infty} w(r)\d r < 0
$$
for $w_0'$ sufficiently small.
\end{thm}

\section{A Necessary Geometric Condition}
The fact that the surface has a reflection symmetry through the $xy$-plane actually provides some geometric conditions on the surface.  Suppose we have a function $z(r)$ with $w=z'=z_r$ satisfying the equation (\ref{VEqn1}), we must check that the surface obtained by revolving the function $-z(r)$ is also stationary to the functional $\HelF$.    For this purpose, we may parametrize the surface by an even function $r(z)$ with $z \in [-a,a]$, namely,
$$
\left\{\begin{aligned}
\bX(\theta,z) &= (r(z)\cos\theta, r(z)\sin\theta,z), \qquad\quad
(\theta,z) \in [0,2\pi]\times [-a,a], \\
r(z) &=r(-z), \quad r(0) = r_\infty.
\end{aligned}\right.
$$
In this section, we will use $z_r$ or $w_r$ to stand for derivatives wrt $r$, i.e., $z'$ and $w'$ previously in order to distinguish from derivatives wrt $z$ denoted by another subscript.

The Helfrich functional for $r(z)$ and $-r(z)$ near $r_\infty$ can be written as,
$$
\frac{\HelF}{2\pi} = \lim_{\varepsilon\to 0}  \int_{-a}^{-\varepsilon} + \int_{\varepsilon}^{a}\, \left(\,\left[(2H+c_0)^2\right] r\sqrt{1+r_z^2} + \frac{p}{2} r^2 \,\right)\,\d z.
$$

We then study the critical function of the functional.  Assume that the variation $\delta(r(z)) = \varphi(z)$, then the above is written as
\begin{gather*}
\frac{\delta(\HelF)}{2\pi} = \lim_{\varepsilon\to 0} \int \left( F_2 \varphi_{zz} + F_1 \varphi_z + F_0 \varphi \right)\d z ;\\
\intertext{where (with similar calculations as those in the
appendix)}
\begin{aligned}
F_2 &= \frac{2r r_{zz}}{(r_z^2+1)^{5/2}} - \frac{2}{(r_z^2+1)^{3/2}} + \frac{2c_0 r}{r_z^2+1} ;\\
F_1 &= \frac{-5rr_z r_{zz}^2}{(r_z^2+1)^{7/2}} - \frac{r_z}{r(r_z^2+1)^{3/2}} + \frac{(c_0^2+\lambda)r r_z} {(r_z^2+1)^{1/2}} + \frac{6 r_z r_{zz}}{(r_z^2+1)^{5/2}} - \frac{4c_0 r r_z r_{zz}}{(r_z^2+1)^2}; \\
F_0 &= \frac{r_{zz}^2}{(r_z^2+1)^{5/2}} - \frac{1}{r^2 (r_z^2+1)^{1/2}} + \frac{2c_0 r_{zz}}{r_z^2+1} + (c_0^2+\lambda)\sqrt{r_z^2+1} + pr.
\end{aligned}
\end{gather*}
Since $F_2$ and $F_0$ are even and $F_1$ is odd, after integrating by parts, we have
$$
\frac{\delta(\HelF)}{2\pi} = \int_{-a}^a \left[ (F_2)_{zz} - (F_1)_z + F_0 \right] \varphi \d z - \lim_{\varepsilon\to 0} \left[ -(F_2)_z(\varepsilon) + F_1(\varepsilon) \right]\left[ \varphi_z(\varepsilon) - \varphi_z(-\varepsilon)\right].
$$
The function $r(z)$ here is the inverse function of the function $z(r)$ in the previous sections.  Thus, $r_z$ and $r_{zz}$ above can be written in terms of $w$ by chain rule.  Since $w(r)=z'(r)$ is a solution to the equation~(\ref{VEqn1}), we obtain
$$
(F_2)_{zz} - (F_1)_z + F_0 = 0.
$$
Furthermore, we may denote $\eta(r) \DEF -(F_2)_z + F_1$.  Then, in terms of $w$,
$$
\eta = \frac{rw_r^2}{(1+w^2)^{5/2}} - \frac{w^2}{r\sqrt{1+w^2}} - 2c_0w - (c_0^2+\lambda) r \sqrt{1+w^2} + \frac{pr^2w}{2}.
$$
Consequently, $\delta(\HelF) = 0$ if and only if $\eta(r)$ is bounded in a neighborhood of $r_\infty$.
Note that
$$
\eta_z = \frac{\eta_r}{w} = \frac{-w_r^2}{w(1+w^2)^{5/2}} + \frac{w}{r^2\sqrt{1+w^2}} - \frac{2c_0w_r}{w(1+w^2)} - (c_0^2+\lambda)\frac{\sqrt{1+w^2}}{w} + pr
$$
is bounded according to the analysis in the previous section, specifically by proposition~\ref{prop-wprimeord}.  Therefore, $\eta(r)$ is also bounded and
hence, we have established the following theorem.
\begin{thm}
Let $z=z(r)\geq 0$ be a function on $[0,r_\infty]$ with $z(r)=0$ if and
only if $r=r_\infty$.  If $w=z'$ is a solution to the
equation~{\em(\ref{VEqn1})}, then the surface obtained by revolving the
curves $z(r)$ and $-z(r)$ is stationary to the functional $\HelF$.
\end{thm}
Furthermore, according to this boundary condition, by considering $\displaystyle\lim_{r\to r_\infty}\frac{\eta(r)}{\sqrt{1+w^2}} = 0$, there is a requirement on the geometry of the surface.
\begin{prop}
Let $\Sigma$ be an axisymmetric stationary surface of $\HelF$.  If $r_\infty$ is the radius of the circle of intersection of $\Sigma$ and the reflection plane, then the Gaussian curvature of every point along this circle is given by
$\displaystyle
K(r_\infty)^2 = \frac{-1}{r_\infty}Q\left( \frac{-1}{r_\infty} \right).
$
\end{prop}

\section{Appendix: Derivation of the Main Equation}
Let the surface of revolution $\Sigma$ be parametrized by
${\mathbf X}(r,\theta) = (r\cos\theta,
r\sin\theta, z(r))$ for $(r,\theta) \in [0,r_\infty] \times
[0,2\pi]$.  Then, with $w(r) = z'(r)$, the mean curvature $H$ is given by
$$
H = \frac{1}{2(1+w^2)} \left[ w' + \frac{1}{r}w + \frac{1}{r}w^3
\right].
$$
The area element is $2\pi r\sqrt{1+w^2}\d r$.  With the assumption that
$z(r_\infty) = 0$, after integrating by
parts, the volume enclosed by $\Sigma$ is
$\displaystyle -\pi \int_0^{r_\infty} r^2 w(r) \d r.$
Therefore, to find the critical point of $\HelF$, we may consider
the variation on
$$
\frac{1}{2\pi}\HelF = \int_0^{r_\infty} \left[ (2H+c_0)^2 + \lambda
\right] r\sqrt{1+w^2}\d r - \frac{p}{2} \int_0^{r_\infty} r^2w \d
r.
$$
Let $\delta$ be the variational operator and $\delta(w) = \varphi$
where $\varphi(r)$ is a smooth function with compact support.
Then $\delta\HelF/2\pi = {\mathcal I}_1 + {\mathcal I}_2 + {\mathcal
I}_3$ where
\begin{align*}
{\mathcal I}_1 &= \int_0^{r_\infty} 2(2H+c_0)\, \delta(H)\,
r\sqrt{1+w^2} \d r; \\
{\mathcal I}_2 &= \int_0^{r_\infty} \left[ (2H+c_0)^2 + \lambda
\right] r\delta(\sqrt{1+w^2})\d r \\
&= \int_0^{r_\infty} \left[ (2H+c_0)^2 + \lambda
\right] r\, \frac{rw\varphi}{\sqrt{1+w^2}} \d r; \\
{\mathcal I}_3 &= \frac{-p}{2}\int_0^{r_\infty} r^2 \delta(w)\d r
=  \frac{-p}{2}\int_0^{r_\infty} r^2 \varphi\d r.
\end{align*}
We also have
\begin{align*}
\delta(2H) &= \frac{1}{(1+w^2)^{5/2}} \left[ (1+w^2)\varphi' +
\frac{1}{r}(1+w^2)\varphi - 3ww'\varphi \right]; \\
(2H)' &= \frac{1}{(1+w^2)^{5/2}} \left[ (1+w^2)w'' - 3w{w'}^2 +
\frac{1}{r}(1+w^2)w' - \frac{1}{r^2}w(1+w^2)^2 \right].
\end{align*}
Therefore,
\begin{align*}
{\mathcal I}_1 &= \int_0^{r_\infty} \frac{2r(2H+c_0)}{(1+w^2)^2}
\left[ (1+w^2)\varphi' +
\frac{1}{r}(1+w^2)\varphi - 3ww'\varphi \right] \\
&= -\int_0^{r_\infty} \pdiff{}{r}\left[ \frac{2r(2H+c_0)}{(1+w^2)}
\right]\varphi \d r + \int_0^{r_\infty}
\frac{2r(2H+c_0)}{(1+w^2)}\varphi \d r \\
& \qquad\qquad {}- \int_0^{r_\infty}
\frac{6rww'(2H+c_0)}{(1+w^2)^2}\varphi \d r \\
&= \int_0^{r_\infty} \left[ \frac{-2r(2H)'}{(1+w^2)} -
\frac{2rww'(2H+c_0)}{(1+w^2)^2} \right]\varphi \d r \\
&= \int_0^{r_\infty} \left[ \frac{-2rw''}{(1+w^2)^{5/2}} +
\frac{6rw{w'}^2}{(1+w^2)^{7/2}} - \frac{2w'}{(1+w^2)^{5/2}} +
\frac{2w}{r(1+w^2)^{3/2}} \right. \\
&\qquad\qquad \left. {} - \frac{2c_0rww'}{(1+w^2)^2} -
\frac{2rw{w'}^2}{(1+w^2)^{7/2}} - \frac{2w^2w'}{(1+w^2)^{5/2}}
\right]\varphi \d r.
\end{align*}
The quantity ${\mathcal I}_2$ is expanded into the following.
\begin{align*}
{\mathcal I}_2 &= \int_0^{r_\infty} 
\frac{\left[(2H+c_0)^2+\lambda\right]rw}{(1+w^2)^{1/2}}\varphi\d r
\\
&= \int_0^{r_\infty} \left[ \frac{rw}{(1+w^2)^{7/2}}\left( w' +
\frac{1}{r}w(1+w^2) \right)^2 \right. \\
&\qquad\qquad\qquad \left. {} + \frac{2c_0rw}{(1+w^2)^2} \left( w' +
\frac{1}{r}w(1+w^2) \right) + \frac{(c_0^2+\lambda)rw}{\sqrt{1+w^2}}
\right]\varphi\d r \\
&=  \int_0^{r_\infty} \left[ \frac{rw{w'}^2}{(1+w^2)^{7/2}} +
\frac{2w^2w'}{(1+w^2)^{5/2}} + \frac{w^3}{r(1+w^2)^{3/2}} \right.
\\
&\qquad\qquad \left. {} + \frac{2c_0rww'}{(1+w^2)^2} +
\frac{2c_0w^2}{(1+w^2)} + \frac{(c_0^2+\lambda)rw}{\sqrt{1+w^2}}
\right]\varphi \d r.
\end{align*}
Summing the above quantities, we have
\begin{align*}
\frac{1}{2\pi}\delta(\HelF) &= {\mathcal I}_1 + {\mathcal I}_2 +
{\mathcal I}_3 \\
&= \int_0^{r_\infty} \left[ \frac{-2rw''}{(1+w^2)^{5/2}} +
\frac{5rw{w'}^2}{(1+w^2)^{7/2}} - \frac{2w'}{(1+w^2)^{5/2}} +
\frac{2w+w^3}{r(1+w^2)^{3/2}} \right. \\
&\qquad\qquad\quad \left. {} + \frac{2c_0w^2}{(1+w^2)} +
\frac{(c_0^2+\lambda)rw}{(1+w^2)^{1/2}} - \frac{pr^2}{2}
\right]\varphi\d r.
\end{align*}
Therefore, the variational equation is given by $(\text{integrand})
+ \text{constant} = 0$.  With the assumption that $w(0)=0$, taking $r\to
0$, one gets
$$
-2w'(0) + \lim_{r\to 0}\frac{2w}{r} + \text{constant} = 0.
$$
This concludes the derivation of the variational equation.


\end{document}

